\documentclass[12pt]{article}
\usepackage{latexsym}

\newcommand {\be}{\begin{equation}}
\newcommand {\ee}{\end{equation}}
\newcommand{\bea}{\begin{eqnarray}}
\newcommand{\eea}{\end{eqnarray}}
\newcommand{\ba}{\begin{array}}
\newcommand{\ea}{\end{array}}
\newcommand{\beq}{\begin{eqnarray*}}
\newcommand{\eeq}{\end{eqnarray*}}

\def\beq{\begin{equation}}
\def\eeq{\end{equation}}
\def\beqn{\begin{eqnarray}}
\def\eeqn{\end{eqnarray}}
\newcommand{\Uqso}{{U'_q({\rm so}_n)}}

\begin{document}

\centerline
{\bf
A $q$-ANALOGUE OF THE WIGNER--ECKART  }

\centerline
{\bf
THEOREM FOR THE NONSTANDARD  }

\centerline
{\bf $q$-DEFORMED ALGEBRA $U'_q({\rm so}_n)$ }

\medskip
\centerline
{N.Z. Iorgov }

\medskip
\centerline
{Bogolyubov Institute for Theoretical Physics}

\centerline
{Metrolohichna Str., 14-b, Kyiv-143, Ukraine, UA-03143}

\centerline
{e-mail: mmtpitp@bitp.kiev.ua}

\vspace{1cm}

\centerline{\bf Abstract}

\noindent
{\small
The tensor product of vector and arbitrary representations
of the nonstandard $q$-deformation $U'_q({\rm so}_n)$
of the universal enveloping algebra
$U({\rm so}_n)$ of Lie algebra ${\rm so}_n$ is defined.
The Clebsch--Gordan coefficients of tensor
product of vector and arbitrary classical or nonclassical type
representations
of $q$-algebra $U'_q({\rm so}_n)$ are found in an explicit form.
The Wigner--Eckart theorem for vector operators is proved.
}
\bigskip

\noindent
{\bf 1. Introduction}

\medskip

\noindent
For the last fifteen years, much attention of mathematicians and
mathematical physicists is attracted to the subject of quantum
algebras and quantum groups. Besides the standard deformation
of Lie algebras proposed by Drinfeld [\ref{Dr}] and Jimbo [\ref{J}], other
(nonstandard) deformations are also under consideration. This paper deals
with the deformation $\Uqso$ of universal enveloping algebra
$U({\rm so}_n)$ proposed by Gavrilik and Klimyk [\ref{GK}].
Let us mention that the algebra $U'_q({\rm so}_3)$
appeared earlier in the paper [\ref{F}].

As a matter of interest, the algebras $\Uqso$ arose naturally
as auxiliary algebras in deriving the algebra of observables in
2+1 quantum gravity with 2D space of genus $g$, so that $n$ depends on $g$,
$n=2g+2$ [\ref{NR},\ref{G},\ref{ChF}].

As shown in [\ref{GI}],
the algebra $\Uqso$ admits a $q$-analogue of Gel'fand--Tsetlin
formalism for construction of finite-dimensional irreducible
representations.
Since the algebra $\Uqso$ is not a Hopf algebra, there is no a natural
way to introduce the notion of tensor product of representations.
But, as shown in [\ref{N},\ref{NUW},\ref{IK:QR}],
the algebra $\Uqso$ is a subalgebra
in Drinfeld--Jimbo Hopf algebra $U_q({\rm sl}_n)$.
Moreover, it is possible to show that the algebra $\Uqso$ is a
$U_q({\rm sl}_n)$-comodule algebra such that the coaction coincides
with the comultiplication in $U_q({\rm sl}_n)$ if one embeds
$\Uqso$ into $U_q({\rm sl}_n)$. This comodule structure can be used
to introduce the tensor product of vector and arbitrary representations
$T$ of $\Uqso$ (it will be denoted by $T^\otimes$), see [\ref{CGC_JPA}].

We describe the decomposition of $T^\otimes$ into
irreducible subrepresentations and write down the corresponding
Clebsch--Gordan coefficients in the case when $T$ is irreducible
finite-dimensional representation of the classical or nonclassical type.
The decomposition of $T^\otimes$ in the case
of classical type representations
has the same form as in the case of Lie algebra ${\rm so}_n$
and the corresponding Clebsch--Gordan coefficients are $q$-deformation
of their classical analogues [\ref{G_SO},\ref{KME}].

It is well-known that Wigner--Eckart theorem for the tensor
operators with respect to Lie algebra ${\rm so}_n$ (and,
especially, ${\rm so}_3$) is very important in physics.
In this paper, we give a $q$--analogue of such
theorem for the case of vector operators.

Everywhere below we suppose that $q$ is not a root of unity.

\bigskip
\noindent{\bf 2. The $q$-deformed algebra $\Uqso$ and quantum algebra
$U_q({\rm sl}_n)$}

\medskip

\noindent
According to [\ref{GK}], the nonstandard $q$-deformation $\Uqso$
of the Lie  algebra ${\rm so}_n$ is
given as a complex associative algebra with ${n-1}$
generating elements $I_{21}$, $I_{32}, \ldots$, $I_{n,n-1}$
obeying the defining relations
\beq
\begin{array}{l}
 I_{j,j-1}^2I_{j-1,j-2} + I_{j-1,j-2}I_{j,j-1}^2 -
[2] \ I_{j,j-1}I_{j-1,j-2}I_{j,j-1} = -I_{j-1,j-2}, \\[2mm]
 I_{j-1,j-2}^2I_{j,j-1} + I_{j,j-1}I_{j-1,j-2}^2 -
[2] \ I_{j-1,j-2}I_{j,j-1}I_{j-1,j-2} = -I_{j,j-1},  \\[2mm]
 [I_{i,i-1},I_{j,j-1}] =0  \qquad {\rm if} \quad \mid {i-j}\mid >1,
\end{array}  \label{f1}
\eeq
where $q+q^{-1}\equiv [2]$, $q\in {\bf C}$, $q\ne 0,\pm 1$.
%%%%%%%%%%%%%%%%%%%%%%%%%%%%%%%%%%%%%%%%%%%%%%%%
It is useful to introduce the generators
\beq                                                     \label{Ipm}
I^{\pm}_{k,l}\equiv [I_{l+1,l} , I^{\pm}_{k,l+1}]_{q^{\pm 1}},
\qquad\qquad k > l+1, \ \  1\leq k,l \leq n,
\eeq
where $[X,Y]_{q^{\pm 1}}\equiv q^{\pm 1/2} X Y - q^{\mp 1/2} Y X$
and $I^+_{k+1,k}\equiv I^-_{k+1,k}\equiv I_{k+1,k}$.
%%%%%%%%
If $q\to 1$ (`classical' limit), the set of relations
(\ref{f1}) reduce to those of $U({\rm so}_n)$.

The algebra $\Uqso$ can be embedded into
quantum algebra $U_q({\rm sl}_n)$, which is defined
[\ref{Dr},\ref{J},\ref{KS}] as a complex associative algebra
with the generating elements $e_i,f_i,k_i,k_i^{-1}$,
$i=1,2,\ldots,n-1$, and defining relations
\[
k_i k_i^{-1}=k_i^{-1} k_i=1, \ \ \ k_i k_j=k_j k_i,  \ \ \
k_ie_jk_i^{-1}=q^{a_{ij}}e_j,\ \ \ k_if_jk_i^{-1}=q^{-a_{ij}}f_j,
\]
\[
[e_i,e_j]=[f_i,f_j]=0,\ \ \ \vert i-j\vert >1,\qquad
[e_i,f_j]=\delta _{ij}{k_i-k_i^{-1}\over q-q^{-1}},
\]
\[
e_i^2e_{i\pm 1}-(q+q^{-1})e_ie_{i\pm 1}e_i+e_{i\pm 1}e^2_i=0, \ \ \
f_i^2f_{i\pm 1}-(q+q^{-1})f_if_{i\pm 1}f_i+f_{i\pm 1}f^2_i=0,
\]
where $a_{ii}=2$, $a_{i,i\pm 1}=-1$ and $a_{ij}=0$ for
$\vert i-j\vert >1$.
It is shown in [\ref{N},\ref{NUW}], that the elements
$\tilde{I}_{i+1,i}=f_i-q^{-1}k_i e_i$,
$i=1,2,\ldots,n-1$,
satisfy the relations (\ref{f1}) and
define a homomorphism  $\Uqso \to U_q({\rm sl}_n)$.
Moreover, it is proved in [\ref{IK:QR}] that this
homomorphism is an embedding, that is, we may consider $\Uqso$
as a subalgebra in $U_q({\rm sl}_n)$.

The quantum algebra $U_q({\rm sl}_n)$ possesses the Hopf structure.
Comultiplication on generators of this algebra can be defined as
\[
 \Delta(e_i)=e_i\otimes k_i^{-1}+1\otimes e_i, \qquad
\Delta(f_i)=f_i\otimes 1+ k_i\otimes f_i, \qquad
\Delta(k_i)=k_i\otimes k_i.
\]
Therefore, we obtain the coideal property of $\Uqso$ embedded into
$U_q({\rm sl}_n)$:
\[
\Delta (\tilde{I}_{i+1,i})=\tilde{I}_{i+1,i}\otimes 1
+ k_i\otimes \tilde{I}_{i+1,i}.
\]

\noindent
{\bf Proposition 1.} {\it The algebra $\Uqso$ is a
$U_q({\rm sl}_n)$-comodule algebra with the coaction
$\phi(I_{i+1,i})=\tilde{I}_{i+1,i}\otimes 1
+ k_i\otimes I_{i+1,i}$. If one embeds $\Uqso$ into $U_q({\rm sl}_n)$,
coaction $\phi$ reduces to comultiplication $\Delta$ of $U_q({\rm sl}_n)$.}

\noindent
{\bf Proof.} This proposition can be verified by direct calculation.
\hfill $\Box$

In particular, Proposition~1 claims that $\phi$ is a homomorphism
from $\Uqso$ into $U_q({\rm sl}_n)\otimes \Uqso$.
This comodule structure can be used
to introduce the tensor product of vector and arbitrary representations
$T$ of $\Uqso$.

Let $T$ be a representation of $\Uqso$ on the linear space $\cal V$
with the basis $\{v_\alpha\}$ and ${\cal V}_{\bf 1}$ be the
$n$-dimensional linear space with the basis $\{v_k\}$,
$k=1,2,\ldots,n$, and ${\cal V}^\otimes\equiv{\cal V}_{\bf 1}\otimes{\cal
V}$.

\noindent
{\bf Proposition 2.} {\it The map $T^\otimes$ from $\Uqso$
to $\mathop{\rm End} {\cal V}^\otimes$
given by the formulas
\begin{equation}        \label{tact1}
T^\otimes(I_{j,j-1})\,(v_{j-1} \otimes v_\alpha)= q\ v_{j-1} \otimes
T(I_{j,j-1})v_\alpha - q^{1/2}\ v_{j} \otimes v_\alpha,
\end{equation}
\begin{equation}        \label{tact2}
T^\otimes(I_{j,j-1})\,(v_{j} \otimes v_\alpha)= q^{-1}\ v_{j} \otimes
T(I_{j,j-1})v_\alpha + q^{-1/2}\ v_{j-1} \otimes v_\alpha,
\end{equation}
\begin{equation}        \label{tact3}
T^\otimes(I_{j,j-1})\,(v_k \otimes v_\alpha)= v_k \otimes
T(I_{j,j-1})v_\alpha, \qquad
\mbox{$j\ne k$,\  $j-1\ne k$}
\end{equation}
defines a representation of $\Uqso$ on the space
${\cal V}^\otimes$. }

\noindent
{\bf Proof.}
Let us define representation ${\cal T}_{\bf 1}$ of
$U_q({\rm sl}_n)$ on the space ${\cal V}_{\bf 1}$
by the formulas
\begin{equation}   \label{sl_act}
\begin{array}{l}
{\cal T}_{\bf 1}(e_i)\,v_k=-q^{-1/2} \delta_{i+1,k} v_{k-1},
\quad
{\cal T}_{\bf 1}(f_i)\,v_k=-q^{1/2} \delta_{i,k} v_{k+1},\\
{\cal T}_{\bf 1}(k_i)\,v_k=q^{\delta_{i,k}-\delta_{i+1,k}} v_k.
\end{array}
\end{equation}
It is easy to verify that this representation is the vector
representation (that is, representation with the highest
weight $(1,0,\ldots,0)$).
The action formulas (\ref{sl_act}) imply
\begin{equation} \label{so_sl_act}
{\cal T}_{\bf 1}(\tilde{I}_{i+1,i})\,v_k=-q^{1/2}\delta_{i,k} v_{k+1}
+q^{-1/2}\delta_{i+1,k} v_{k-1}.
\end{equation}
This representation of $\Uqso$
is equivalent to the classical type
repre\-sen\-tation $T_{{\bf m}_n}$
with ${\bf m}_n=(1,0,\ldots,0)$, that is, the vector representation
(see next section). Hence, similarly to the classical case,
the restriction of the vector representation of $U_q({\rm sl}_n)$
onto $\Uqso$ is the vector representation of $\Uqso$.
This proposition immediately follows from Proposition~1 and
formula (\ref{so_sl_act}), if one takes
$T^\otimes=({\cal T}_1\otimes T)\circ \phi$.
\hfill $\Box$

In the case when $T$ is the trivial representation of $\Uqso$ given by
formulas
$T(a)=0$, $a\in \Uqso$, $a\ne 1$,
Proposition~2 gives us a representation on the
space ${\cal V}_{\bf 1}\sim{\cal V}^\otimes$. We denote this representation
by $T_{\bf 1}$.
\[
T_{\bf 1}(I_{j,j-1})\,v_k=-q^{1/2}\delta_{k,j-1} v_{j}
+q^{-1/2}\delta_{k,j} v_{j-1}.
\]
The representations $T_{\bf 1}$ and $T_{{\bf m}_n}$,
${\bf m}_n=(1,0,\ldots,0)$ (see next section), are equivalent.

In the limit $q\to 1$, Proposition~2 defines the representation which is
the tensor product of the vector and some arbitrary representation of the
Lie algebra ${\rm so}_n$.
On the base of these two arguments, we shall also use the notion
$T^\otimes\equiv T_{\bf 1}\otimes T$.

\bigskip
\noindent{\bf 3. Finite dimensional classical type
representations of $\Uqso$}

\medskip

\noindent
In this section we describe
(in the framework of Gel'fand--Tsetlin formalism)
irreducible finite-dimensional
representation of the algebra $U'_q({\rm so}_{n})$,
which are $q$-deformations of the finite-dimensional irreducible
representations of the Lie algebra ${\rm so}_n$.
They are given by sets ${\bf m}_{n}$
consisting of $\lfloor n/2 \rfloor$ numbers
$m_{1,n}$, $m_{2,n},\ldots$, $m_{\lfloor n/2\rfloor ,n}$
(here $\lfloor n/2 \rfloor$ denotes integral part
of $n/2$) which are all integral or all half-integral and
satisfy the dominance conditions
\begin{equation} \label{DC}
\begin{array}{l}
m_{1,2p+1}{\ge}m_{2,2p+1}{\ge}...{\ge}
m_{p,2p+1}{\ge 0},  \\
m_{1,2p}{\ge}m_{2,2p}{\ge}...{\ge}m_{p-1,2p}{\ge}|m_{p,2p}|
\end{array}
\end{equation}
for $n=2p+1$ and $n=2p$, respectively.
These representations are denoted by $T_{{\bf m}_n}$.
For a basis in a representation space ${\cal V}_{{\bf m}_n}$
we take the $q$-analogue of
Gel'fand--Tsetlin basis which is obtained by successive reduction of
the representation $T_{{\bf m}_n}$ to the subalgebras
$U'_q({\rm so}_{n-1})$, $U'_q({\rm so}_{n-2})$, $\cdots$, $U'_q({\rm
so}_3)$,
$U'_q({\rm so}_2)\equiv U({\rm so}_2)$.
As in the classical case, its elements are labelled by Gel'fand--Tsetlin
tableaux
\begin{equation}\label{GT}
  \{\xi_{n} \}
\equiv \{ {\bf m}_{n},\xi_{n-1}\}\equiv \{{\bf m}_{n} ,
{\bf m}_{n-1} ,\xi_{n-2}\}\equiv \cdots
\equiv \{ {\bf m}_{n}, {\bf m}_{n-1}, \ldots,  {\bf m}_{2}\},
\end{equation}
where the components of ${\bf m}_{k}$ and ${\bf m}_{k-1}$ satisfy the
``betweenness'' conditions
\[
\begin{array}{l}
 m_{1,2p+1}\ge m_{1,2p}\ge m_{2,2p+1} \ge m_{2,2p} \ge \ldots
\ge m_{p,2p+1} \ge m_{p,2p} \ge -m_{p,2p+1}  , \\[2mm]
 m_{1,2p}\ge m_{1,2p-1}\ge m_{2,2p} \ge m_{2,2p-1} \ge \ldots
\ge m_{p-1,2p-1} \ge \vert m_{p,2p} \vert .
\end{array}
\]
The basis element defined by tableau $\{\xi_{n} \}$ is denoted as
$|\xi_{n} \rangle$.
We suppose that the representation space ${\cal V}_{{\bf m}_n}$
is a Hilbert space
and vectors $|\xi_{n} \rangle$ are orthonormal.
It is convenient to introduce the so-called $l$-coordinates
\begin{equation}
l_{j,2p+1}=m_{j,2p+1}+p-j+1,  \qquad
                      l_{j,2p}=m_{j,2p}+p-j         \label{lcoord}
\end{equation}
for the numbers $m_{i,k}$.
The operator $T_{{\bf m}_n}(I_{2p+1,2p})$ of the representation
$T_{{\bf m}_n}$ of $U'_q({\rm so}_{n})$ acts upon Gel'fand--Tsetlin
basis elements, labeled by (\ref{GT}), as
\begin{equation} \label{Act1}
 T_{{\bf m}_n}(I_{2p+1,2p})
| \xi_n\rangle =
\sum^p_{j=1} A^j_{2p}(\xi_n)
            \vert (\xi_n)^{+j}_{2p}\rangle -
\sum^p_{j=1} A^j_{2p}((\xi_n)^{-j}_{2p})
|(\xi_n)^{-j}_{2p}\rangle
\end{equation}
and the operator $T_{{\bf m}_n}(I_{2p,2p-1})$ of the representation
$T_{{\bf m}_n}$ acts as
\begin{equation} \label{Act2}
\begin{array}{l}
 T_{{\bf m}_n}(I_{2p,2p-1})\vert \xi_n\rangle=
\sum^{p-1}_{j=1} B^j_{2p-1}(\xi_n)
\vert (\xi_n)^{+j}_{2p-1} \rangle
 \\[2mm]
-\sum^{p-1}_{j=1} B^j_{2p-1}((\xi_n)^{-j}_{2p-1})
\vert (\xi_n)^{-j}_{2p-1}\rangle
+ {\rm i}\, C_{2p-1}(\xi_n) \vert \xi_n \rangle , \\[3mm]
 T_{{\bf m}_n}(I_{21})\vert \xi_n\rangle=
{\rm i}\, [l_{12}] \vert \xi_n\rangle.
\end{array}
\end{equation}
In these formulas, $(\xi_n)^{\pm j}_{k}$ means the tableau (\ref{GT})
in which $j$-th component $m_{j,k}$ in ${\bf m}_k$ is replaced
by $m_{j,k}\pm 1$. The coefficients
$A^j_{2p},  $ $B^j_{2p-1},$ $C_{2p-1}$
in (\ref{Act1}) and (\ref{Act2}) are given
by the expressions
\begin{equation} \label{ME_A}
A^j_{2p}(\xi_n) =
\left(
\frac{[l_{j,2p}] [l_{j,2p}+1]}
{[2 l_{j,2p}] [2 l_{j,2p}+2]}\right)^{\frac12} {\hat A}^j_2p\ ,
\end{equation}
\begin{equation} \label{ME_B}
B^j_{2p-1}(\xi_n)=\frac{{\hat B}^j_{2p-1}(\xi_n)}
{ [l_{j,2p-1}] \bigl([2 l_{j,2p-1}+1] [2 l_{j,2p-1}-1]\bigr)^\frac12 } ,
\end{equation}
\[
{\hat A}^j_{2p}=\left(
\frac{ \prod_{i=1}^p
[l_{i,2p+1}+l_{j,2p}][l_{i,2p+1}-l_{j,2p}-1] }
 {  \prod_{i\ne j}^p
[l_{i,2p}+l_{j,2p}] [l_{i,2p}-l_{j,2p}] }  \right.
\]
\begin{equation} \label{ME_AH}
 \times \left. \frac{   \prod_{i=1}^{p-1}
[l_{i,2p-1}+l_{j,2p}][l_{i,2p-1}-l_{j,2p}-1] }
 {  \prod_{i\ne j}^p
[l_{i,2p}+l_{j,2p}+1] [l_{i,2p}-l_{j,2p}-1]  }
\right)^{\frac12}
\end{equation}
and
\[
{\hat B}^j_{2p-1}(\xi_n)=
\left(
\frac { \prod_{i=1}^p[l_{i,2p}+l_{j,2p-1}][l_{i,2p}-l_{j,2p-1}] }
{     \prod_{i\ne j}^{p-1}
 [l_{i,2p-1}+l_{j,2p-1}] [l_{i,2p-1}-l_{j,2p-1}] }
\right.
\]
\begin{equation} \label{ME_BH}
\times \left.
\frac { \prod_{i=1}^{p-1}
[l_{i,2p-2}+l_{j,2p-1}][l_{i,2p-2}-l_{j,2p-1}] }
  {\prod_{i\ne j}^{p-1}
[l_{i,2p-1}+l_{j,2p-1}-1] [l_{i,2p-1}-l_{j,2p-1} -1]  }
\right)^{\frac12} ,
\end{equation}

\begin{equation} \label{ME_C}
 C_{2p-1}(\xi_n) =\frac{ \prod_{i=1}^p [ l_{i,2p} ]
\prod_{i=1}^{p-1} [ l_{i,2p-2} ]}
{\prod_{i=1}^{p-1} [l_{i,2p-1}] [l_{i,2p-1} - 1] } ,
\end{equation}
where numbers in square brackets mean $q$-numbers defined by
$[a]:= (q^a-q^{-a})/(q-q^{-1})$.

\bigskip
\noindent{\bf 4. Finite dimensional nonclassical type
representations of $\Uqso$}

\medskip

\noindent
The representations of the previous section are called representations
of the classical type, because at $q\to 1$ the operators
$T_{{\bf m}_n}(I_{j,j-1})$ turn into the corresponding operators
$T_{{\bf m}_n}(I_{j,j-1})$ for irreducible finite dimensional
representations with highest weights ${\bf m}_n$ of the Lie algebra
${\rm so}_n$.

The algebra $U'_q({\rm so}_n)$ also has irreducible finite dimensional
representations $T$ of nonclassical type, that is, such that the operators
$T(I_{j,j-1})$ have no classical limit $q\to 1$.
They are given (see [\ref{IK:NC}])
by sets $\epsilon := (\epsilon _2,\epsilon _3,\cdots ,
\epsilon _n)$, $\epsilon _i=\pm 1$, and by sets
${\bf m}_{n}$ consisting of $\lfloor{n/2}\rfloor $ {\bf half-integral}
numbers $m_{1,n}$, $m_{2,n}, \ldots,$ $m_{\lfloor n/2\rfloor ,n}$
(here $\lfloor {n/2}\rfloor$ denotes integral part of $n/2$)
that satisfy the dominance conditions
\begin{equation} \label{DCN}
m_{1,n}\ge m_{2,n}\ge ... \ge
m_{\lfloor n/2\rfloor,n}\ge 1/2.
\end{equation}
These representations are denoted by $T_{\epsilon,{\bf m}_n}$.

For a basis in the representation space
${\tilde{\cal V}}_{{\bf m}_n}$ we use the analogue of the
basis of the previous section. Its elements are
labeled by tableaux
\begin{equation}\label{GTN}
  \{\xi_{n} \}
\equiv \{ {\bf m}_{n},\xi_{n-1}\}\equiv \{{\bf m}_{n} ,
{\bf m}_{n-1} ,\xi_{n-2}\}\equiv \cdots
\equiv \{ {\bf m}_{n}, {\bf m}_{n-1}, \ldots,  {\bf m}_{2}\},
\end{equation}
where the components of ${\bf m}_{k}$ and ${\bf m}_{k-1}$ satisfy the
``betweenness'' conditions
\[
\begin{array}{l}
m_{1,2p+1}\ge m_{1,2p}\ge m_{2,2p+1} \ge m_{2,2p} \ge ...
\ge m_{p,2p+1} \ge m_{p,2p} \ge 1/2  ,\\
m_{1,2p}\ge m_{1,2p-1}\ge m_{2,2p} \ge m_{2,2p-1} \ge ...
\ge m_{p-1,2p-1} \ge m_{p,2p}  .
\end{array}
\]
The basis element defined by tableau $\{\xi_{n} \}$ is denoted
as $\vert \xi_{n} \rangle $.
We suppose that the representation space $\tilde{{\cal V}}_{{\bf m}_n}$
is a Hilbert space
and vectors $|\xi_{n} \rangle$ are orthonormal.
It is convenient to introduce the $l$-coordinates as in (\ref{lcoord}).

The operator $T_{\epsilon,{\bf m}_n}(I_{2p+1,2p})$ of the representation
$T_{\epsilon,{\bf m}_n}$ of $\Uqso$ acts upon
basis elements, labeled by (\ref{GTN}), by the formula
\[
T_{\epsilon,{\bf m}_n}(I_{2p+1,2p})
| \xi_n\rangle =
\delta_{m_{p,2p},1/2}\, \frac{\epsilon _{2p+1}}{q^{1/2}-q^{-1/2}} D_{2p}
(\xi _n) | \xi_n\rangle
\]
\begin{equation} \label{Act1e}
+\sum^p_{j=1} {\tilde A}^j_{2p}(\xi_n)
            \vert (\xi_n)^{+j}_{2p}\rangle -
\sum^p_{j=1} {\tilde A}^j_{2p}((\xi_n)^{-j}_{2p})
|(\xi_n)^{-j}_{2p}\rangle ,
\end{equation}
where the summation in the last sum must be from 1 to $p-1$ if
$m_{p,2p}=1/2$,
and the operator $T_{{\bf m}_n}(I_{2p,2p-1})$ of the representation
$T_{{\bf m}_n}$ acts as
\[
T_{\epsilon ,{\bf m}_n}(I_{2p,2p-1})\vert \xi_n\rangle=
\sum^{p-1}_{j=1} {\tilde B}^j_{2p-1}(\xi_n)
\vert (\xi_n)^{+j}_{2p-1} \rangle
\]
\begin{equation} \label{Act2e}
-\sum^{p-1}_{j=1}{\tilde B}^j_{2p-1}((\xi_n)^{-j}_{2p-1})
\vert (\xi_n)^{-j}_{2p-1}\rangle
+ \epsilon _{2p} {\tilde C}_{2p-1}(\xi_n)
\vert \xi_n \rangle ,
\end{equation}
\[
T_{\epsilon ,{\bf m}_n}(I_{21})\vert \xi_n\rangle=\epsilon_2[l_{12}]_+
|\xi_n\rangle,
\]
where $[a]_+:=(q^a+q^{-a})/(q-q^{-1})$.
In these formulas, $(\xi_n)^{\pm j}_{k}$ means the tableau (\ref{GTN})
in which $j$-th component $m_{j,k}$ in ${\bf m}_k$ is replaced
by $m_{j,k}\pm 1.$ Matrix elements
${\tilde A}^j_{2p}$ and ${\tilde B}^j_{2p-1}$ are defined using
formulas (\ref{ME_AH}) and (\ref{ME_BH}):
\[
{\tilde A}^j_{2p}(\xi_n)=\frac{{\hat A}^j_{2p}(\xi_n)}
{\bigl((q^{l_{j,2p}}-q^{-l_{j,2p}})(q^{l_{j,2p}+1}-q^{-l_{j,2p}-1})\bigr)^{\
frac12}}
\]
\[
{\tilde B}^j_{2p-1}(\xi_n)=\frac{{\hat B}^j_{2p-1}(\xi_n)}
{[l_{j,2p-1}]_+ \bigl([2l_{j,2p-1}+1][2l_{j,2p-1}-1]\bigr)^{\frac12}},
\]
\[
{\tilde C}_{2p-1}(\xi_n) = \frac {
\prod_{s=1}^p [ l_{s,2p} ]_+
\prod_{s=1}^{p-1} [ l_{s,2p-2} ]_+}{
   \prod_{s=1}^{p-1} [l_{s,2p-1}]_+ [l_{s,2p-1} - 1]_+ } .
\]
\[
D_{2p} (\xi _n)=
\frac{\prod_{i=1}^p
[l_{i,2p+1}-\frac 12 ] \prod_{i=1}^{p-1} [l_{i,2p-1}-\frac 12 ] }
{\prod_{i=1}^{p-1}
[l_{i,2p}+\frac 12 ] [l_{i,2p}-\frac 12 ] } .
\]

\bigskip
\noindent{\bf 5. Decomposition of representations
$T_{\bf 1}\otimes T_{{\bf m}_3}$ of the
algebra $U'_q({\rm so}_3)$}

\medskip

\noindent
In this and in the next sections,
we consider the decomposition of
representations $T^\otimes\equiv T_{\bf 1}\otimes T_{{\bf m}_n}$
into irreducible constituents of the algebra $U'_q({\rm so}_n)$.
In this section, we restrict ourselves to the case $n=2,3$.

First, we consider the case of the algebra
$U'_q({\rm so}_2)\equiv U({\rm so}_2)$. This
algebra has representations
$T_m$, $m\equiv m_{12}$, $m\in\frac12{\bf Z}$, of the {\it classical type}
acting on one-dimensional spaces with basis vectors $|m\rangle$,
and $T_m(I_{21})|m\rangle={\rm i} [m] |m\rangle$.
Then
\[
T^\otimes(I_{21})(v_1\otimes|m\rangle)={\rm i} q [m] v_1\otimes|m\rangle
- q^{1/2} v_2\otimes|m\rangle,
\]
\[
T^\otimes(I_{21})(v_2\otimes|m\rangle)={\rm i} q^{-1} [m]
v_2\otimes|m\rangle
+ q^{-1/2} v_1\otimes|m\rangle.
\]
This representation is reducible. We introduce the vectors
\beq
v_\pm^{(m)}=\mp {\rm i} q^{-1/2\pm m} v_1 + v_2.
\label{vpm}
\eeq
Then the vectors $|m\pm 1\rangle^\otimes:= v_\pm^{(m)} \otimes |m\rangle$
are
eigenvectors of $T^\otimes(I_{21})$:\\
$T^\otimes(I_{21})|m\pm 1\rangle^\otimes={\rm i} [m\pm 1] |m\pm
1\rangle^\otimes$.
This fact can be easy verified by direct calculation using
the definition of $q$-numbers. Thus, we have decomposition
$T^\otimes\equiv T_{\bf 1} \otimes T_m= T_{m+1} \oplus T_{m-1}$.

Now, we consider the case of the algebra $U'_q({\rm so}_3)$. This
algebra has representations $T_l$, ${\bf m}_3\equiv (m_{13})\equiv (l)$,
$l\in \{0,$ $1/2,$ $1,$ $3/2,\ldots\}$, of the {\it classical type} acting
on the
spaces ${\cal V}_l$ with the basis vectors
$|l,m\rangle$, ($m\equiv m_{12}$), $m=-l,-l+1,\ldots,l$:
\[
 T_l(I_{21})|l,m\rangle={\rm i} [m] |l,m\rangle,\qquad
T_l(I_{32})|l,m\rangle=A_{l,m}|l,m+1\rangle-A_{l,m-1}|l,m-1\rangle,
\]
where $A_{l,m}=d_m ([l-m][l+m+1])^{1/2}$,
$d_m=\bigl([m][m+1]/([2m][2m+2])\bigr)^{1/2}$.
Let us consider the vectors
\begin{equation}                                   \label{Sv_so3}
|l',m\rangle^\otimes:= \alpha^{(l')}_{l,m} v_+^{(m-1)}\otimes |l,m-1\rangle
+
\beta^{(l')}_{l,m} v_3\otimes |l,m\rangle +
\gamma^{(l')}_{l,m} v_-^{(m+1)}\otimes |l,m+1\rangle,
\end{equation}
where $m=-l',-l'+1,\ldots,l',$ and
\[
l'=l+1,l,l-1\ \ \ \mbox{if\ \ }l\ge 1;\qquad
l'=3/2,1/2\ \ \ \mbox{if\ \ }l=1/2;
\qquad l'=1\ \ \ \mbox{if\ \ }l=0.
\]
The vectors $v_\pm^{(m)}$ in (\ref{Sv_so3}) are defined in (\ref{vpm}) and
\[\alpha^{(l+1)}_{l,m}= q^{l-m+1/2}d_{m-1}([l+m][l+m+1])^{1/2},\]
\[\beta^{(l+1)}_{l,m}= ([l-m+1][l+m+1])^{1/2},\]
\[\gamma^{(l+1)}_{l,m}= -q^{l+m+1/2}d_{m}([l-m][l-m+1])^{1/2},\]
\[\alpha^{(l)}_{l,m}= -q^{-m-1/2}d_{m-1}([l+m][l-m+1])^{1/2},\]
\[\beta^{(l)}_{l,m}=[m],\]
\[\gamma^{(l)}_{l,m}= -q^{m-1/2}d_{m}([l-m][l+m+1])^{1/2},\]
\[\alpha^{(l-1)}_{l,m}=-q^{-l-m-1/2}d_{m-1}([l-m][l-m+1])^{1/2},\]
\[\beta^{(l-1)}_{l,m}= ([l-m][l+m])^{1/2},\]
\[\gamma^{(l-1)}_{l,m}= q^{-l+m-1/2}d_{m}([l+m][l+m+1])^{1/2}.\]
From the case of $U'_q({\rm so}_2)$, it is easy to see that
$T^\otimes(I_{21})|l',m\rangle^\otimes={\rm i} [m] |l',m\rangle^\otimes$.
One can show by direct calculation that
$T^\otimes(I_{32})|l',m\rangle^\otimes=A_{l',m}|l',m+1\rangle^\otimes
-A_{l',m-1}|l',m-1\rangle^\otimes$.
It means that the vectors $|l',m\rangle^\otimes$ at fixed $l'$ span a
subspace in
${\cal V}^\otimes$, which is invariant and irreducible under the action
of $T^\otimes(a)$, $a\in U'_q({\rm so}_3)$. The corresponding
subrepresentation is equivalent to $T_{l'}$.
Comparing the dimensions of $T_{l'}$ with dimension of
$T^\otimes$, we conclude that
$T^\otimes=T_{l+1}\oplus T_l\oplus T_{l-1}$, if $l\ge 1$;
$T^\otimes=T_{3/2}\oplus T_{1/2}$, if $l=1/2$; $T^\otimes=T_1$, if $l=0$.
Let us remind that $T_{l}\equiv T_{{\bf m}_{3}}$, $m_{13}\equiv l$.
The numbers $\alpha^{(l')}_{l,m}$, $\beta^{(l')}_{l,m}$ and
$\gamma^{(l')}_{l,m}$ are Clebsch--Gordan coefficients of these
decompositions.

\bigskip
\noindent{\bf 6. Decomposition of $T_{\bf 1}\otimes T_{{\bf m}_n}$ of
the algebra $U'_q({\rm so}_n)$, $n\ge 4$}

\medskip

\noindent
In this section, we consider the decomposition of the
representations $T^\otimes \equiv T_{\bf 1}\otimes T_{{\bf m}_n}$
of algebra $U'_q({\rm so}_n)$, $n\ge 4$, into irreducible constituents.
All the results of this section are obtained in [\ref{CGC_JPA}].
As shown there, this decomposition has the form
\beq
T^\otimes =
\bigoplus_{{\bf m}'_n\in {\cal S}({\bf m}_n)}  T_{{\bf m}'_n},
\label{tn}
\eeq
where
\beq
{\cal S}({\bf m}_{2p+1})=\bigcup_{j=1}^p \{{\bf m}^{+j}_{2p+1}\}
\cup \bigcup_{j=1}^p \{{\bf m}^{-j}_{2p+1}\}
\cup  \{{\bf m}^{\vphantom{+j}}_{2p+1}\},
                   \label{t2p1}
\eeq
\beq
{\cal S}({\bf m}_{2p})=\bigcup_{j=1}^p \{{\bf m}^{+j}_{2p}\}
\cup \bigcup_{j=1}^p \{{\bf m}^{-j}_{2p}\}.
                    \label{t2p}
\eeq
By ${\bf m}_n^{\pm j}$ we mean here the set ${\bf m}_n$ with
$m_{j,n}$ replaced by $m_{j,n}\pm 1$, respectively.
If some ${\bf m}_n^{\pm j}$ is not dominant (\ref{DC}),
then the corresponding ${\bf m}_n^{\pm j}$
must be omitted. If $m_{p,2p+1}=0$ then ${\bf m}_{2p+1}$ in
right-hand side of (\ref{t2p1}) also must be omitted.
For decomposition (\ref{tn}) of the representation
$T^\otimes$, there correspond the decomposition
of carrier space:
\begin{equation} \label{vn}
{\cal V}^\otimes \equiv
{\cal V}_{\bf 1}\otimes {\cal V}_{{\bf m}_n}=
\bigoplus_{{\bf m}'_n\in{\cal S}({\bf m}_n)}
{\cal V}_{{\bf m}'_n}.
\end{equation}
In order to give this decomposition in an explicit form,
we change the basis $\{v_k\otimes |\xi_n\rangle\}$, $k=1,2,\ldots,n$,
in ${\cal V}^\otimes$ to
$\{v_k\otimes |\xi_n\rangle\}$, $k=+,-,3,\ldots,n$, by replacing (for
every fixed $\{\xi_n\}=
\{{\bf m}_n,{\bf m}_{n-1}, \ldots, {\bf m}_3, {\bf m}_2\}$)
two basis vectors $v_1\otimes |\xi_n\rangle$ and
$v_2\otimes |\xi_n\rangle$ by
$v_+^{(m_{12})}\otimes |\xi_n\rangle$ and
$v_-^{(m_{12})}\otimes |\xi_n\rangle$
(see (\ref{vpm})).
From now on, we shall omit the index $(m_{12})$ in the notion
of the basis vectors $v_\pm^{(m_{12})}\otimes |\xi_n\rangle$,
supposing that it is equal to $m_{12}$-component of the corresponding
Gel'fand--Tsetlin tableaux $\{\xi_n\}$.

We introduce the vectors
(where $\{\xi'_n\}=\{{\bf m}'_n,{\bf m}'_{n-1}, \ldots, {\bf m}'_3,{\bf
m}'_2\}$)
\beq
 |{\bf m}'_n,\xi_{n-1}\rangle^\otimes:=
\sum_k \sum_{|{\bf m}_n,\xi'_{n-1}\rangle\in {\cal V}_{{\bf m}_n}}
\bigl(k, ({\bf m}_n,\xi'_{n-1})|({\bf m}'_n,\xi_{n-1})\bigr)
\,v_k\otimes |{\bf m}_n,\xi'_{n-1}\rangle
\label{SvCGC}
\eeq
in the space ${\cal V}^\otimes$, where $k$ runs over the set
$+,-,3,\ldots,n$, and coefficients\\
$\bigl(k, ({\bf m}_n,\xi'_{n-1})|({\bf m}'_n,\xi_{n-1})\bigr)$
are Clebsch--Gordan coefficients (CGC's).
Now we define these CGC's in an explicit form.

We put $\bigl(k, ({\bf m}_n,\xi'_{n-1})|
({\bf m}'_n,\xi_{n-1})\bigr)=0$ if one of the conditions
\[
\begin{array}{l}
 {\rm 1)\ } {\bf m}'_n\not\in{\cal S}({\bf m}_n), \\
{\rm 2)\ } {\bf m}_s\not\in{\cal S}({\bf m}'_s), s=n-1,\ldots,k,
\ \ k\ge 3,\\
 {\rm 3)\ } {\bf m}_s\not\in{\cal S}({\bf m}'_s),
s=n-1,\ldots,3,\ \ k=+,-,\\
{\rm 4)\ } \xi'_{k-1}\ne\xi_{k-1}, k=3,4,\ldots,n, \\
 {\rm 5)\ } m_{12} \ne m'_{12}+1, k=+,\\
{\rm 6)\ } m_{12} \ne m'_{12}-1, k=-.
\end{array}
\]
is fulfilled.
The nonzero CGC for $k=n$ are:
\beq                   \label{SCGC1}
\begin{array}{l}
\bigl(2p+1, ({\bf m}_{2p+1},\xi_{2p})|
({\bf m}^{+j}_{2p+1},\xi_{2p})\bigr){=}
\Bigl(\prod_{r=1}^p
[l_{j,2p+1}+l_{r,2p}][l_{j,2p+1}-l_{r,2p}]\Bigr)^{\frac12},\\[3mm]
\bigl(2p+1, ({\bf m}_{2p+1},\xi_{2p})|
({\bf m}_{2p+1},\xi_{2p})\bigr){=}
\prod_{r=1}^p [l_{r,2p}],\\[1mm]
\bigl(2p+1, ({\bf m}_{2p+1},\xi_{2p})|
({\bf m}^{-j}_{2p+1},\xi_{2p})\bigr){=}
\Bigl(\prod_{r=1}^p
[l_{j,2p+1}{+}l_{r,2p}{-}1][l_{j,2p+1}{-}l_{r,2p}{-}1]\Bigr)^{\frac12},
\end{array}
\eeq
\beq                          \label{SCGC0}
\begin{array}{l}
\bigl(2p, ({\bf m}_{2p},\xi_{2p-1})|
({\bf m}^{+j}_{2p},\xi_{2p-1})\bigr){=}
\Bigl(\prod_{r=1}^{p-1}
[l_{j,2p}+l_{r,2p-1}][l_{j,2p}-l_{r,2p-1}+1]\Bigr)^{\frac12},\\[2mm]
\bigl(2p, ({\bf m}_{2p},\xi_{2p-1})|
({\bf m}^{-j}_{2p},\xi_{2p-1})\bigr){=}
\Bigl(\prod_{r=1}^{p-1}
[l_{j,2p}+l_{r,2p-1}-1][l_{j,2p}-l_{r,2p-1}]\Bigr)^{\frac12}.
\end{array}
\eeq
(They are defined up to normalization, that is, multiplication of these
CGC's
by some constants will not spoil the following results.)

All the other CGC's can be found from just presented
as follows:
\[
\bigl(k, \xi_n|\xi'_n\bigr)=q^{k-n}
\frac{\langle{\bf m}_{n+1},\xi_n|T_{{\bf m}_{n+1}}
(I_{n+1,k}^-)|{\bf m}_{n+1},\xi'_n\rangle}
{\langle{\bf m}_{n+1},{\bf m}_n,\xi_{n-1}|T_{{\bf m}_{n+1}}
(I_{n+1,n})|{\bf m}_{n+1},
{\bf m}'_n,\xi_{n-1}\rangle}
\]
\beq
\times
\bigl(n, ({\bf m}_n,\xi_{n-1})|({\bf m}'_n,\xi_{n-1})\bigr),
                                    \label{CGCk}
\eeq
where the generators $I_{n+1,k}^-$ are defined in (\ref{Ipm}).
If $k=+$ or $k=-$ in the left-hand side of (\ref{CGCk}),
one must put $k=2$ in right-hand side.
The set ${\bf m}_{n+1}$ must be chosen to give non-zero denominator
in right-hand side of (\ref{CGCk}).
Note that if $(n,({\bf m}_n,\xi_{n-1})|({\bf m}'_n,\xi_{n-1}))\ne 0$,
one  can always do such a choice, moreover, the resulting
CGC will not depend on this particular choice.
In the case $n=3$ we reobtain the CGC's
for the algebra $U'_q({\rm so}_3)$ (see section 5).

As shown in [\ref{CGC_JPA}], the defined CGC's have the
{\it factorization} property. This fact (in complete analogy
with the classical case, see [\ref{G_SO},\ref{KME}]) gives a possibility
to present arbitrary CGC for the algebra $U'_q({\rm so}_{n})$
as a product of {\it scalar factors}.

\noindent {\bf Theorem 1.}
{\it
The formulas for the action of the
operators $T^\otimes(I_{k+1,k})$, $k=1,2,\ldots, n-1$, on the vectors
$|{\bf m}'_n,\xi_{n-1}\rangle^\otimes$
defined by (\ref{SvCGC}) with
CGC's defined by (\ref{SCGC1})--(\ref{CGCk}),
coincide with the corresponding formulas (\ref{Act1})--(\ref{Act2})
for the action of the operators
$T_{{\bf m}'_n}(I_{k+1,k})$ on the GT basis
vectors $|{\bf m}'_n,\xi_{n-1}\rangle$.
We have the decomposition (\ref{tn}).
}

\newpage

\noindent{\bf 7. Decomposition of representations
$T_{\bf 1}\otimes T_{\epsilon, {\bf m}_3}$ of the
algebra $U'_q({\rm so}_3)$}

\medskip

\noindent
In this and in the next section,
we consider the decomposition of
representations $T^\otimes\equiv T_{\bf 1}\otimes T_{\epsilon,{\bf m}_n}$
into irreducible constituents of the algebra $U'_q({\rm so}_n)$.
In this section, we restrict ourselves to the case $n=2,3$.

First, we consider the case of the algebra
$U'_q({\rm so}_2)\equiv U({\rm so}_2)$. This
algebra has representations $T_{\epsilon_2,m}$, $\epsilon_2=\pm 1$,
$m\equiv m_{12}$, $m\in \{1/2, 3/2, \ldots\}$,
acting on one-dimensional spaces with basis vectors $|m\rangle$,
and $T_{\epsilon_2, m}(I_{21})|m\rangle= \epsilon_2 [m]_+ |m\rangle$.
Then the representation $T^\otimes\equiv T_{\bf 1}\otimes T_{\epsilon_2, m}$
is two-dimensional and reducible. We introduce the vectors
\beq
v_\pm^{(\epsilon_2,m)}=-\epsilon_2 q^{-1/2\pm m} v_1 + v_2.
\label{vpme}
\eeq
Then the vectors $|m\pm 1\rangle^\otimes:= v_\pm^{(\epsilon_2,m)}
\otimes |m\rangle$ are
eigenvectors of $T^\otimes(I_{21})$:
$T^\otimes(I_{21})|m\pm 1\rangle^\otimes=\epsilon_2 [m\pm 1]_+
|m\pm 1\rangle^\otimes$.
This fact can be easy verified by direct calculation using
the definition of $q$-numbers. Thus, we have decomposition
$T^\otimes\equiv T_{\bf 1} \otimes T_{\epsilon,m}=
T_{\epsilon_2,m+1} \oplus T_{\epsilon_2,m-1}$, if $m\ge 3/2$, and
$T^\otimes\equiv T_{\bf 1} \otimes T_{\epsilon_2,1/2}=
T_{\epsilon_2,3/2} \oplus T_{\epsilon_2,1/2}$.

Now, we consider the case of the algebra $U'_q({\rm so}_3)$. This
algebra has four classes of representations of {\it nonclassical type}
$T_{\epsilon,l}$, $\epsilon=\{\epsilon_2,\epsilon_3\}$,
$\epsilon_i\in\{\pm 1\}$,
 ${\bf m}_3\equiv (m_{13})\equiv (l)$,
$l\in \{1/2,$ $3/2,$ $5/2,  \ldots\}$, acting on the
spaces ${\cal V}_l$ with the basis vectors
$|l,m\rangle$, ($m\equiv m_{12}$), $m=1/2,3/2,\ldots,l$:
\[
 T_{\epsilon,l}(I_{21})|l,m\rangle=\epsilon_2 [m]_+ |l,m\rangle,
\]
\[
T_{\epsilon,l}(I_{32})|l,m\rangle=\tilde{A}_{l,m}|l,m+1\rangle-
\tilde{A}_{l,m-1}|l,m-1\rangle,
\qquad \mbox{if $m\ge 3/2$},
\]
\[
T_{\epsilon,l}(I_{32})|l,1/2\rangle=\tilde{A}_{l,1/2}|l,3/2\rangle
+\epsilon_3 [1/2]_+ [l+1/2] |l,1/2\rangle,
\]
where $\tilde {A}_{l,m}=\tilde{d}_m ([l-m][l+m+1])^{1/2}$,
$\tilde{d}_m=\bigl((q^m-q^{-m})(q^{m+1}-q^{-m-1})\bigr)^{-1/2}$.
Let us consider the vectors
\begin{equation}                                   \label{Sv_so3e}
|l',m\rangle^\otimes:=
\tilde{\alpha}^{(l')}_{l,m} v_+^{(\epsilon_2,m-1)}\otimes |l,m-1\rangle +
\tilde{\beta}^{(l')}_{l,m} v_3\otimes |l,m\rangle +
\tilde{\gamma}^{(l')}_{l,m} v_-^{(\epsilon_2,m+1)}\otimes |l,m+1\rangle,
\end{equation}
where $m=3/2,5/2, \ldots,l',$ and
\[
l'=l+1,l,l-1\ \ \ \mbox{if\ \ }l\ge 3/2;\qquad
l'=3/2,1/2\ \ \ \mbox{if\ \ }l=1/2.
\]
If $m=1/2$, we should replace $|l,-1/2\rangle$ by  $|l,1/2\rangle$
in right-hand side of (\ref{Sv_so3e}).

The vectors $v_\pm^{(\epsilon_2,m)}$ in (\ref{Sv_so3e}) are
defined in (\ref{vpme}) and
\[\tilde{\alpha}^{(l+1)}_{l,m}=
q^{l-m+1/2}\tilde{d}_{m-1}([l+m][l+m+1])^{1/2},
\qquad m\ne 1/2\]
\[\tilde{\beta}^{(l+1)}_{l,m}= ([l-m+1][l+m+1])^{1/2},\]
\[\tilde{\gamma}^{(l+1)}_{l,m}= -q^{l+m+1/2}\tilde{d}_{m}([l-m][l-m+1])^{1/2
},\]
\[\tilde{\alpha}^{(l)}_{l,m}= q^{-m-1/2}\tilde{d}_{m-1}([l+m][l-m+1])^{1/2},
\qquad m\ne 1/2\]
\[\tilde{\beta}^{(l)}_{l,m}=[m]_+,\]
\[\tilde{\gamma}^{(l)}_{l,m}= -q^{m-1/2}\tilde{d}_{m}([l-m][l+m+1])^{1/2},\]
\[\tilde{\alpha}^{(l-1)}_{l,m}=-q^{-l-m-1/2}\tilde{d}_{m-1}([l-m][l-m+1])^{1
/2},
\qquad m\ne 1/2\]
\[\tilde{\beta}^{(l-1)}_{l,m}= ([l-m][l+m])^{1/2},\]
\[\tilde{\gamma}^{(l-1)}_{l,m}=
q^{-l+m-1/2}\tilde{d}_{m}([l+m][l+m+1])^{1/2},\]
\[\tilde{\alpha}^{(l+1)}_{l,1/2}= -q^{l}[1/2]_+\epsilon_3
([l+1/2][l+3/2])^{1/2},\]
\[\tilde{\alpha}^{(l+1)}_{l,1/2}= -q^{-1}[1/2]_+\epsilon_3 [l+1/2],\]
\[\tilde{\alpha}^{(l-1)}_{l,1/2}= q^{-l-1}[1/2]_+\epsilon_3
([l-1/2][l+1/2])^{1/2}.\]
From the case of $U'_q({\rm so}_2)$, it is easy to see that
$T^\otimes(I_{21})|l',m\rangle^\otimes=$ $\epsilon_2 [m]_+
|l',m\rangle^\otimes$.
One can show by direct calculation that the operator $T^\otimes(I_{32})$
acts on the set of vectors $|l',m\rangle^\otimes$ at some fixed $l'$
as operator $T_{\epsilon,l'}(I_{32})$ acts on the Gel'fand--Tsetlin
basis vectors $|l',m\rangle$.
It means that the vectors $|l',m\rangle^\otimes$ at fixed $l'$ span a
subspace in
${\cal V}^\otimes$, which is invariant and irreducible under the action
of $T^\otimes(a)$, $a\in U'_q({\rm so}_3)$. The corresponding
subrepresentation is equivalent to $T_{\epsilon,l'}$.
Comparing the dimensions of $T_{\epsilon,l'}$ with dimension of
$T^\otimes$, we conclude that
$T^\otimes=T_{\epsilon, l+1}\oplus T_{\epsilon, l}\oplus
T_{\epsilon,l-1}$, if $l\ge 3/2$;
$T^\otimes=T_{\epsilon, 3/2}\oplus T_{\epsilon, 1/2}$, if $l=1/2$.
Let us remind that $T_{\epsilon, l}\equiv T_{\epsilon, {\bf m}_{3}}$,
$m_{13}\equiv l$.
The numbers $\tilde{\alpha}^{(l')}_{l,m}$, $\tilde{\beta}^{(l')}_{l,m}$ and
$\tilde{\gamma}^{(l')}_{l,m}$ are Clebsch--Gordan coefficients of these
decompositions.

\bigskip
\noindent{\bf 8. Decomposition of $T_{\bf 1}\otimes T_{\epsilon,{\bf m}_n}$
of
the algebra $U'_q({\rm so}_n)$, $n\ge 4$}

\medskip

\noindent
In this section, we describe the decomposition of the
representations $T^\otimes \equiv T_{\bf 1}\otimes T_{\epsilon, {\bf m}_n}$
of algebra $U'_q({\rm so}_n)$, $n\ge 4$, into irreducible constituents.
This decomposition has the form
\beq
T^\otimes =
\bigoplus_{{\bf m}'_n\in {\cal S}({\bf m}_n)}  T_{\epsilon, {\bf m}'_n},
\label{tne}
\eeq
where
\beq
{\cal S}({\bf m}_{2p+1})=\bigcup_{j=1}^p \{{\bf m}^{+j}_{2p+1}\}
\cup \bigcup_{j=1}^p \{{\bf m}^{-j}_{2p+1}\}
\cup  \{{\bf m}^{\vphantom{+j}}_{2p+1}\},
                   \label{t2p1e}
\eeq
\beq
{\cal S}({\bf m}_{2p})=\bigcup_{j=1}^p \{{\bf m}^{+j}_{2p}\}
\cup \bigcup_{j=1}^p \{{\bf m}^{-j}_{2p}\}.
                    \label{t2pe}
\eeq
By ${\bf m}_n^{\pm j}$ we mean here the set ${\bf m}_n$ with
$m_{j,n}$ replaced by $m_{j,n}\pm 1$, respectively.
If $m_{p,2p}=1/2$, the element ${\bf m}^{-p}_{2p}$
in right-hand side of (\ref{t2pe}) must be replaced by ${\bf m}_{2p}$.
If some ${\bf m}_n^{\pm j}$ is not dominant (\ref{DCN}),
then the corresponding ${\bf m}_n^{\pm j}$
must be omitted; in particular, if $m_{p,2p+1}=1/2$, the element
${\bf m}^{-p}_{2p+1}$ must be omitted.
Note, that the representation $T_{\bf 1}\otimes T_{\epsilon, {\bf m}_n}$
decomposes into irreducible nonclassical type representations
with the same set $\epsilon=(\epsilon_2,\epsilon_3,\ldots)$.
For decomposition (\ref{tne}) of the representation
$T^\otimes$, there correspond the decomposition
of carrier space:
\begin{equation} \label{vne}
{\cal V}^\otimes \equiv
{\cal V}_{\bf 1}\otimes {\tilde{\cal V}}_{\epsilon, {\bf m}_n}=
\bigoplus_{{\bf m}'_n\in{\cal S}({\bf m}_n)}
{\tilde{\cal V}}_{\epsilon, {\bf m}'_n}.
\end{equation}

In order to give this decomposition in an explicit form,
we change the basis $\{v_k\otimes |\xi_n\rangle\}$, $k=1,2,\ldots,n$,
in ${\cal V}^\otimes$ to
$\{v_k\otimes |\xi_n\rangle\}$, $k=+,-,3,\ldots,n$, by replacing (for
every fixed $\{\xi_n\}=
\{{\bf m}_n,{\bf m}_{n-1}, \ldots, {\bf m}_3, {\bf m}_2\}$)
two basis vectors $v_1\otimes |\xi_n\rangle$ and
$v_2\otimes |\xi_n\rangle$ by
$v_+^{(\epsilon_2, m_{12})}\otimes |\xi_n\rangle$ and
$v_-^{(\epsilon_2, m_{12})}\otimes |\xi_n\rangle$
(see (\ref{vpme})).
From now on, we shall omit the index $(\epsilon_2, m_{12})$ in the notion
of the basis vectors $v_\pm^{(\epsilon_2, m_{12})}\otimes |\xi_n\rangle$,
supposing that it contains $m_{12}$-component of the corresponding
Gel'fand--Tsetlin tableaux $\{\xi_n\}$.
% and $\epsilon_2$ is fixed to be
%the first component of $\epsilon$.

We introduce the vectors
(where $\{\xi'_n\}=\{{\bf m}'_n,{\bf m}'_{n-1}, \ldots, {\bf m}'_3,{\bf
m}'_2\}$)
\beq
 |{\bf m}'_n,\xi_{n-1}\rangle^\otimes:=
\sum_k \sum_{|{\bf m}_n,\xi'_{n-1}\rangle\in {\cal V}_{{\bf m}_n}}
\bigl(k, ({\bf m}_n,\xi'_{n-1})|({\bf m}'_n,\xi_{n-1}); \epsilon\bigr)
\,v_k\otimes |{\bf m}_n,\xi'_{n-1}\rangle
\label{SvCGCe}
\eeq
in the space ${\cal V}^\otimes$, where $k$ runs over the set
$+,-,3,\ldots,n$, and coefficients\\
$\bigl(k, ({\bf m}_n,\xi'_{n-1})|({\bf m}'_n,\xi_{n-1}); \epsilon \bigr)$
are Clebsch--Gordan coefficients (CGC's).
Now we define these CGC's in
an explicit form.

We put $\bigl(k, ({\bf m}_n,\xi'_{n-1})|({\bf
m}'_n,\xi_{n-1}); \epsilon\bigr)=0$ if one of the conditions
\[
\begin{array}{l}
 {\rm 1)\ } {\bf m}'_n\not\in{\cal S}({\bf m}_n), \\
 {\rm 2)\ } {\bf m}_s\not\in{\cal S}({\bf m}'_s), s=n-1,\ldots,k,
\ \ k\ge 3, \\
 {\rm 3)\ } {\bf m}_s\not\in{\cal S}({\bf m}'_s),
s=n-1,\ldots,3,\ \ k=+,-,\\
{\rm 4)\ } \xi'_{k-1}\ne\xi_{k-1}, k=3,4,\ldots,n, \\
{\rm 5)\ } m_{12} \ne m'_{12}+1, k=+,\\
{\rm 6)\ } m_{12} \ne m'_{12}-1,\ \  k=-,\ \  m'_{12}\ge \frac 32,\\
{\rm 6')\ } m_{12} \ne m'_{12},\ \  k=-,\ \  m'_{12}=\frac 12.
\end{array}
\]
is fulfilled.
The nonzero CGC for $k=n$ are:
\beq                   \label{SCGC1e}
\begin{array}{l}
\bigl(2p{+}1, ({\bf m}_{2p+1},\xi_{2p})|
({\bf m}^{+j}_{2p+1},\xi_{2p});\epsilon\bigr){=}
\Bigl(\prod_{r=1}^p
[l_{j,2p+1}+l_{r,2p}][l_{j,2p+1}-l_{r,2p}]\Bigr)^{\frac12},\\[3mm]
\bigl(2p{+}1, ({\bf m}_{2p+1},\xi_{2p})|
({\bf m}_{2p+1},\xi_{2p});\epsilon\bigr){=}
\prod_{r=1}^p [l_{r,2p}]_+,\\[1mm]
\bigl(2p{+}1, ({\bf m}_{2p+1},\xi_{2p})|
({\bf m}^{-j}_{2p+1},\xi_{2p});\epsilon\bigr){=}
\Bigl(\prod_{r=1}^p
[l_{j,2p+1}{+}l_{r,2p}{-}1][l_{j,2p+1}{-}l_{r,2p}{-}1]\Bigr)^{\frac12},
\end{array}
\eeq
\beq                          \label{SCGC0e}
\begin{array}{l}
\bigl(2p, ({\bf m}_{2p},\xi_{2p-1})|
({\bf m}^{+j}_{2p},\xi_{2p-1});\epsilon\bigr){=}
\Bigl(\prod_{r=1}^{p-1}
[l_{j,2p}+l_{r,2p-1}][l_{j,2p}-l_{r,2p-1}+1]\Bigr)^{\frac12},\\[2mm]
\bigl(2p, ({\bf m}_{2p},\xi_{2p-1})|
({\bf m}^{-j}_{2p},\xi_{2p-1});\epsilon\bigr){=}
\Bigl(\prod_{r=1}^{p-1}
[l_{j,2p}+l_{r,2p-1}-1][l_{j,2p}-l_{r,2p-1}]\Bigr)^{\frac12},\\[2mm]
\bigl(2p, ({\bf m}_{2p},\xi_{2p-1})|
({\bf m}_{2p},\xi_{2p-1});\epsilon\bigr){=}
\prod_{r=1}^{p-1}
[l_{r,2p-1}-\frac 12], \qquad \mbox{if }m_{p,2p}=\frac 12.
\end{array}
\eeq
(They are defined up to normalization, that is, multiplication of these
CGC's
by some constants will not spoil the following results.)

All the other CGC's can be found from just presented
by the following formula:
\[
\bigl(k, \xi_n|\xi'_n;\epsilon\bigr)=q^{k-n}
\frac{\langle{\bf m}_{n+1},\xi_n|T_{\tilde\epsilon, {\bf m}_{n+1}}
(I_{n+1,k}^-)|{\bf m}_{n+1},\xi'_n\rangle}
{\langle{\bf m}_{n+1},{\bf m}_n,\xi_{n-1}|T_{\tilde\epsilon, {\bf m}_{n+1}}
(I_{n+1,n})|{\bf m}_{n+1},
{\bf m}'_n,\xi_{n-1}\rangle}
\]
\beq
\times
\bigl(n, ({\bf m}_n,\xi_{n-1})|({\bf m}'_n,\xi_{n-1});\epsilon\bigr),
                                    \label{CGCke}
\eeq
where the generators $I_{n+1,k}^-$ are defined in (\ref{Ipm}),
$\tilde\epsilon=(\epsilon_2, \epsilon_3, \ldots, \epsilon_n, +1)$.
If $k=+$ or $k=-$ in the left-hand side of (\ref{CGCke}),
one must put $k=2$ in right-hand side.
The set ${\bf m}_{n+1}$ must be chosen to give non-zero denominator
in right-hand side of (\ref{CGCke}).
Note that if $(n,({\bf m}_n,\xi_{n-1})|({\bf m}'_n,\xi_{n-1}))\ne 0$,
one  can always do such a choice, moreover, the resulting
CGC will not depend on this particular choice.
In the case $n=3$ we reobtain the CGC's corresponding to the
nonclassical type representations
for the algebra $U'_q({\rm so}_3)$ (see section 7).
\smallskip

\noindent {\bf Theorem 2.}
{\it
The formulas for the action of the
operators $T^\otimes(I_{k+1,k})$, $k=1,2,\ldots, n-1$, on the vectors
$|{\bf m}'_n,\xi_{n-1}\rangle^\otimes$
defined by (\ref{SvCGCe}) with
CGC's defined by (\ref{SCGC1e})--(\ref{CGCke}),
coincide with the corresponding formulas (\ref{Act1e})--(\ref{Act2e})
for the action of the operators
$T_{\epsilon,{\bf m}'_n}(I_{k+1,k})$ on the GT basis
vectors $|{\bf m}'_n,\xi_{n-1}\rangle$.
We have the decomposition (\ref{tne}).
}

\bigskip
\noindent{\bf 9. The Wigner--Eckart theorem for the vector operators}

\medskip

\noindent
To fix idea, we restrict ourselves to the case when vector operator
acts on the space where direct sum of classical
type representations of $\Uqso$ is realized.

The formula (\ref{SvCGC}) give us the transformation from the
basis $\{v_k\otimes |\xi_{n}\rangle\}$ to the basis
$\{|\xi'_{n}\rangle^\otimes\}$ in the space ${\cal V}^\otimes$.
Because of (\ref{vn}), the transformation matrix is non-degenerate matrix
with matrix elements being CGC's $(k,\xi_n|\xi'_n)$.
Denote the matrix elements of inverse matrix by
$(\xi'_n|k,\xi_n)$ ({\it inverse} CGC's).
Let us find the expression for the vector
$v_n\otimes |{\bf m}_n,\xi_{n-1}\rangle$ from (\ref{SvCGC})
in terms of vectors $|\xi'_{n}\rangle^\otimes$.
Since this vector transforms under
the action of $T^\otimes(a)$, $a\in U'_q({\rm so}_{n-1})$, as the vector
$|\xi_{n-1}\rangle$ under the action of $T_{{\bf m}_{n-1}}(a)$
(see formula (\ref{tact3})), from Schur lemma we have
\begin{equation}
v_n\otimes |{\bf m}_n,\xi_{n-1}\rangle=
\sum_{{\bf m}'_n}
\bigl(({\bf m}'_n,\xi_{n-1})|n, ({\bf m}_n,\xi_{n-1})\bigr)
\,|{\bf m}'_n,\xi_{n-1}\rangle^\otimes,
\label{InvCGC}
\end{equation}
where the coefficients $\bigl(({\bf m}'_n,\xi_{n-1})|n,
({\bf m}_n,\xi_{n-1})\bigr)$ depend only on
${\bf m}'_n$, ${\bf m}_n$, ${\bf m}_{n-1}$.
From (\ref{SvCGC}), it also follows that
${\bf m}'_n\in {\cal S}({\bf m}_n)$. Although these coefficients
are uniquely defined by (\ref{SvCGC})--(\ref{CGCk}), we shall need only
their explicit dependence on ${\bf m}_{n-1}$.

\noindent
{\bf Definition 1.} {\it The set $\{V_k\}$, $k=1,2,\ldots,n$, of
operators on $\cal V$, where a representation $T$ of $\Uqso$ is realized,
such that
\begin{equation}        \label{to1}
[ V_{j-1}, T(I_{j,j-1})]_q = V_{j},\ \ \ [T(I_{j,j-1}), V_{j}]_q = V_{j-1},
\end{equation}
\begin{equation}        \label{to2}
[T(I_{j,j-1}), V_{k}]=0, \qquad
\mbox{if $j\ne k$ and $j-1\ne k$,}
\end{equation}
where $[X,Y]_q{=}q^{1/2}X Y{-}q^{-1/2}Y X$,
is called vector operator of the algebra $\Uqso$. }

It is easy to verify, that the action of operators $T(I_{j,j-1})$
on the vectors $V_k\, v_\alpha$ directly correspond to the action
(\ref{tact1})--(\ref{tact3}) of operators  $T^\otimes(I_{j,j-1})$
on the vectors $v_k \otimes v_\alpha$.

Let $T$ be a direct sum of
irreducible classical type representations of $\Uqso$
with arbitrary multiplicities. Choose Gel'fand--Tsetlin (GT) basis
in ${\cal V}$.
Let us consider an invariant subspace ${\cal V}_{{\bf m}_n,s}$
where subrepresentation equivalent to $T_{{\bf m}_n}$ is
realized.  The number $s$ labels the number of such subspace
if the corresponding multiplicity exceeds $1$.
Combine the vectors $V_k\ |({\bf m}_n,\xi_{n-1}); s\rangle$,
where $\{|({\bf m}_n,\xi_{n-1}); s\rangle\}$ is GT basis of
${\cal V}_{{\bf m}_n,s}$, with CGC as in (\ref{SvCGC}) for some
fixed ${\bf m}'_n\in {\cal S}({\bf m}_n)$. It is possible
two variants. First, all the vectors $|{\bf m}'_n,\xi_{n-1}\rangle^\otimes$
are zero. Second, on the space spanned by the vectors
$|{\bf m}'_n,\xi_{n-1}\rangle^\otimes$, a representation
of $\Uqso$ equivalent to $T_{{\bf m}'_n}$ is realized.
From Schur lemma, it follows that
\begin{equation} \label{WE_a}
|{\bf m}'_n,\xi_{n-1}\rangle^\otimes=\sum_{s'} ({\bf m}'_n,s'\|V\|{\bf
m}_n,s)
\, |{\bf m}'_n,\xi_{n-1}; s'\rangle,
\end{equation}
where $({\bf m}'_n,s'\|V\|{\bf m}_n,s)$ are some coefficients
({\it reduced matrix elements}) depending only on ${\bf m}'_n$, $s'$,
${\bf m}_n$, $s$ and vector operator $\{V_k\}$.
Using the analogue of relation (\ref{InvCGC}) for vector operator
and (\ref{WE_a}) we have
\begin{eqnarray}
V_n |{\bf m}_n,\xi_{n-1}; s\rangle=
\sum_{{\bf m}'_n, s'}
\bigl(({\bf m}'_n,\xi_{n-1})|n, ({\bf m}_n,\xi_{n-1})\bigr)
\nonumber\\
\qquad\times
({\bf m}'_n,s'\|V\|{\bf m}_n,s)
\, |{\bf m}'_n,\xi_{n-1}; s'\rangle.
\label{VnACT}
\end{eqnarray}
As was claimed above, the coefficients
$\bigl(({\bf m}'_n,\xi_{n-1})|n, ({\bf m}_n,\xi_{n-1})\bigr)$
may depend on ${\bf m}_{n-1}$. Since this dependence is identical
for all the possible vector operators in arbitrary spaces, we
choose, for a moment, $\cal V$ to be the space ${\cal V}_{{\bf m}_{n+1}}$
of irreducible representation $T_{{\bf m}_{n+1}}$ of
$U'_q({\rm so}_{n+1})$ for some convenient ${\bf m}_{n+1}$,
and $\{V_k\}\equiv\{T_{{\bf m}_{n+1}}(I^+_{n+1,k})\}$. Extracting
the dependence on ${\bf m}_{n-1}$ from the matrix elements
of $T_{{\bf m}_{n+1}}(I_{n+1,n})$ and comparing it with formulas
(\ref{SCGC1})--(\ref{SCGC0}), we obtain
\[
\bigl(({\bf m}'_n,\xi_{n-1})|n, ({\bf m}_n,\xi_{n-1})\bigr)=
\bigl(n, ({\bf m}_n,\xi_{n-1})|({\bf m}'_n,\xi_{n-1})\bigr)
\lambda_{{\bf m}'_n,{\bf m}_n},
\]
where $\lambda_{{\bf m}'_n,{\bf m}_n}$ are some coefficients
depending on ${\bf m}'_n$ and ${\bf m}_n$ only.
Returning to the formula (\ref{VnACT}) and
denoting $({\bf m}'_n,s'\|V\|{\bf m}_n,s)'=$
$({\bf m}'_n,s'\|V\|{\bf m}_n,s)\times$ $\lambda_{{\bf m}'_n,{\bf m}_n}$
we have
\begin{eqnarray}
V_n |{\bf m}_n,\xi_{n-1}; s\rangle=
\sum_{{\bf m}'_n, s'}
\bigl(n, ({\bf m}_n,\xi_{n-1})|({\bf m}'_n,\xi_{n-1})\bigr)\nonumber\\
\qquad \times ({\bf m}'_n,s'\|V\|{\bf m}_n,s)'
\, |{\bf m}'_n,\xi_{n-1}; s'\rangle.
\end{eqnarray}
Iterating the second formula in (\ref{to1}),
we obtain the action formulas for $\{V_k\}$, $1\leq k<n$.
Thus, we deduce the following $q$-analogue of Wigner--Eckart
theorem.

\noindent {\bf Theorem 4.}
{\it
If $\cal V$ is a Hilbert space and its Gel'fand--Tsetlin basis
\\$\{|{\bf m}_n,\xi_{n-1}; s\rangle\}$ is orthonormal, we have, for
the components of vector operator $\{V_k\}$ on $\cal V$,
the decomposition
\[
\langle {\bf m}'_n,\xi'_{n-1}; s'|
V_k |{\bf m}_n,\xi_{n-1}; s\rangle=
\bigl(({\bf m}'_n,\xi'_{n-1})|k, ({\bf m}_n,\xi_{n-1})\bigr)'
({\bf m}'_n,s'\|V\|{\bf m}_n,s)',
\]
where
\[
\bigl(({\bf m}'_n,\xi'_{n-1})|k, ({\bf m}_n,\xi_{n-1})\bigr)'=
\frac{\langle{\bf m}_{n+1},\xi'_n|T_{{\bf m}_{n+1}}
(I_{n+1,k}^+)|{\bf m}_{n+1},\xi_n\rangle}
{\langle{\bf m}_{n+1},{\bf m}'_n,\xi_{n-1}|T_{{\bf m}_{n+1}}
(I_{n+1,n})|{\bf m}_{n+1},
{\bf m}_n,\xi_{n-1}\rangle}
\]
\[
\qquad \times
\bigl(n, ({\bf m}_n,\xi_{n-1})|({\bf m}'_n,\xi_{n-1})\bigr),
\qquad 1\leq k<n                                    \label{InvCGCk}
\]
(see comments after analogous formula (\ref{CGCk})).
}

Let $T_\epsilon$ be a direct sum of
irreducible nonclassical type representations of $\Uqso$
with arbitrary multiplicities and fixed $\epsilon$ on the
Hilbert space ${\cal V}_\epsilon$.
Choose Gel'fand--Tsetlin basis in ${\cal V}_\epsilon$.
The the space ${\cal V}_\epsilon$ is direct sum of
subspaces ${\cal V}_{\epsilon, {\bf m}_n,s}$
where subrepresentations equivalent to $T_{\epsilon, {\bf m}_n}$ are
realized.  The number $s$ labels the number of such subspace
if the corresponding multiplicity exceeds $1$.
Using argumentation analogous to the case of classical type
representations, we derive the following $q$-analogue of Wigner--Eckart
theorem for the case of nonclassical type representations.

\noindent {\bf Theorem 5.}
{\it
If ${\cal V}_\epsilon$ is a Hilbert space and its Gel'fand--Tsetlin basis\\
$\{|{\bf m}_n,\xi_{n-1}; s\rangle\}$ is orthonormal, we have, for
the components of vector ope\-ra\-tor $\{V_k\}$ on $\cal V_\epsilon$,
the decomposition
\[
\langle {\bf m}'_n,\xi'_{n-1}; s'|
V_k |{\bf m}_n,\xi_{n-1}; s\rangle=
\bigl(({\bf m}'_n,\xi'_{n-1})|k, ({\bf m}_n,\xi_{n-1});\epsilon\bigr)'
\]
\[\times
(\epsilon,{\bf m}'_n,s'\|V\|\epsilon,{\bf m}_n,s)',
\]
where
\[
\bigl(({\bf m}'_n,\xi'_{n-1})|k, ({\bf m}_n,\xi_{n-1});\epsilon\bigr)'=
\]
\[
\frac{\langle{\bf m}_{n+1},\xi'_n|T_{\tilde \epsilon, {\bf m}_{n+1}}
(I_{n+1,k}^+)|{\bf m}_{n+1},\xi_n\rangle}
{\langle{\bf m}_{n+1},{\bf m}'_n,\xi_{n-1}|T_{\tilde \epsilon, {\bf
m}_{n+1}}
(I_{n+1,n})|{\bf m}_{n+1},
{\bf m}_n,\xi_{n-1}\rangle}
\]
\[
\qquad \times
\bigl(n, ({\bf m}_n,\xi_{n-1})|({\bf m}'_n,\xi_{n-1}); \epsilon\bigr),
\qquad 1\leq k<n                                    \label{InvCGCke}
\]
(see comments after analogous formula (\ref{CGCke})).
}

The coefficients $(\epsilon,{\bf m}'_n,s'\|V\|\epsilon, {\bf m}_n,s)'$
are {\it reduced matrix elements} for the vector operator $\{V_k\}$.

If the representation $T$ is a direct sum of classical type
representations and nonclassical type representations with
different $\epsilon$, it is easy to find the matrix elements
for the vector operators. It is sufficient to take into account
the fact that vector operator `acting' on classical type representation
can not give nonclassical type representation, and `acting'
on nonclassical type representation with some set $\epsilon$
can not give nonclassical type representation with other set $\epsilon'$.
Thus, corresponding matrix elements are zero. The non-zero matrix
elements are described by Theorem~4 and Theorem~5.

\medskip
The author is thankful to A.~U.~Klimyk and
A.~M.~Gavrilik for the fruitful discussions.
The research described in this article was made possible in part
by Award No. UP1-2115 of the U.S. Civilian Research and Development
Foundation (CRDF).

\bigskip

\centerline{\bf REFERENCES}

\smallskip

\begin{enumerate}
\item \label{Dr} Drinfeld V. G. Hopf algebra and Yang--Baxter equation //
Sov. Math. Dokl. -- 1985. -- {\bf 32}. -- P.~254--259.

\item\label{J} Jimbo M. A $q$-difference analogue of $U(g)$
and the Yang--Baxter Equation // Lett. Math. Phys. -- 1985. -- {\bf 10}. --
P.~63--69.

\item\label{GK}
Gavrilik A.M., Klimyk A.U. $q$-Deformed orthogonal and
pseudo-ortho\-go\-nal algebras and their representations //
Lett.~Math.~Phys. -- 1991. -- {\bf 21}. -- P.~215--220.

\item\label{F}
Fairlie D.B. Quantum deformation of $SU_q(2)$
// J.~Phys~A. -- 1990. -- {\bf 23}. -- P.~L183--L187.

\item\label{NR} Nelson J., Regge T. $2+1$ Quantum gravity
// Phys. Lett. -- 1991. -- {\bf B272}. --
P.~213--216.

\item\label{G} Gavrilik A.M. The use of quantum algebras in quantum
gravity//
Proc. of Inst. of Math. of NAS Ukraine. -- 2000. -- {\bf 30}. --
P.~304--309.

\item\label{ChF} Chekhov L.O., Fock V.V. Observables in 3D gravity and
geodesic algebras // Czech.~J.~Phys. --- 2000. -- {\bf 50}. --
P.~1201--1208.

\item\label{GI}
Gavrilik A.M., Iorgov N.Z. $q$-Deformed algebras
$U_q(so_n)$ and their representations // Methods Funct. Anal.
Topol. -- 1997. -- {\bf 3}, No4. -- P.~51--63.

\item\label{N}
Noumi M. Macdonald's symmetric polynomials as zonal
spherical functions on some quantum homogeneous spaces // Adv. Math.
-- 1996. -- {\bf 123}. -- P.~16--77.

\item\label{NUW} Noumi M., Umeda T., Wakayama M. Dual pairs,
spherical harmonics and a Capelli identity in quantum group theory //
Compos. Math. -- 1996. -- {\bf 104}. -- P.~227--277.

\item\label{IK:QR} Iorgov N.Z., Klimyk A.U. Nonstandard deformation
$U'_q({\rm so}_n)$ for $q$ a root of unity // Methods Funct. Anal.
Topol. -- 2000. -- {\bf 6}, No3. -- P.~56--71; math.QA/0007105.

\item\label{CGC_JPA} Iorgov N.Z. On tensor products of representations of
the non-standard $q$-deformed algebra $U'_q({\rm so}_n)$ // J.~Phys.~A. --
2001. -- {\bf 34}. -- P.~3095--3108.

\item\label{G_SO} Gavrilik A.M. Clebsch--Gordan coefficients for
the direct product $[m_n]\otimes[1]$ of the group
${\rm SO}(n)$ representations. -- Kiev, 1973 (Preprint ITP-73-155P
of Institute for Theoretical Physics) (in Russian).

\item\label{KME} Klimyk A.U. Matrix Ele\-ments and Clebsch--Gordan
Coeffi\-ci\-ents of
Group Re\-pre\-sentations. -- Kiev: Naukova Dumka, 1979 (in Russian).

\item\label{KS} Klimyk A., Schm\"udgen K. Quantum Groups and
Their Representations. -- Berlin: Springer, 1997.

\item\label{IK:NC} Iorgov N.Z., Klimyk A.U. Nonclassical type
representations of the $q$-deformed algebra
$U'_q({\rm so}_n)$ // Czech.~J.~Phys. -- 2000. -- {\bf 50}. -- P.~85--90.

\end{enumerate}
\end{document}